\documentclass[12pt]{amsart}
\usepackage[nobysame,abbrev,alphabetic]{amsrefs}
\usepackage{amssymb}
\usepackage[all]{xy}
\usepackage{stmaryrd}

\DeclareMathOperator{\cd}{cd}

\DeclareMathOperator{\dirlim}{\varinjlim}

\DeclareMathOperator{\Hom}{Hom}

\DeclareMathOperator{\Supp}{Supp}
\DeclareMathOperator{\trg}{trg}

\newcommand{\Bock}{\mathrm{Bock}}

\newcommand{\nek}{,\ldots,}
\newcommand{\Id}{\mathrm{Id}}

\newcommand{\indec}{\mathrm{indec}}

\newcommand{\inv}{^{-1}}
\newcommand{\isom}{\cong}

\newcommand{\Sh}{\mathrm{Sh}}

\newcommand{\tensor}{\otimes}

\DeclareFontFamily{U}{wncy}{}
\DeclareFontShape{U}{wncy}{m}{n}{<->wncyr10}{}
\DeclareSymbolFont{mcy}{U}{wncy}{m}{n}
\DeclareMathSymbol{\Sha}{\mathord}{mcy}{"58}
\DeclareMathSymbol{\sha}{\mathord}{mcy}{"78}

\newtheorem{thm}{Theorem}[section]
\newtheorem{cor}[thm]{Corollary}
\newtheorem{lem}[thm]{Lemma}
\newtheorem{prop}[thm]{Proposition}
\newtheorem{defin}[thm]{Definition}
\newtheorem{exam}[thm]{Example}

\newtheorem{rem}[thm]{Remark}

\swapnumbers

\numberwithin{equation}{section}

\newcommand{\alp}{\alpha}
\newcommand{\gam}{\gamma}

\newcommand{\eps}{\epsilon}
\newcommand{\lam}{\lambda}
\newcommand{\Lam}{\Lambda}
\newcommand{\sig}{\sigma}

\newcommand{\dbG}{\mathbb{G}}

\newcommand{\dbQ}{\mathbb{Q}}

\newcommand{\dbU}{\mathbb{U}}

\newcommand{\dbZ}{\mathbb{Z}}

\newcommand{\calL}{\mathcal{L}}

\begin{document}

\title[Cohomology for Magnus Formations]{Cohomology and the Combinatorics of Words for Magnus Formations}

\author{ Ido Efrat}
\address{Earl Katz Family Chair in Pure Mathematics\\
Department of Mathematics\\
Ben-Gurion University of the Negev\\
P.O.\ Box 653, Be'er-Sheva 8410501\\
Israel} \email{efrat@bgu.ac.il}

\thanks{This work was supported by the Israel Science Foundation (grant No.\ 569/21).}



\maketitle

\begin{abstract}
For a prime number $p$ and a free pro-$p$ group $G$ on a totally ordered basis $X$, we consider closed normal subgroups $G^\Phi$ of $G$ which are generated by $p$-powers of iterated commutators associated with Lyndon words in the alphabet $X$.
We express the profinite cohomology group $H^2(G/G^\Phi)$ combinatorically,  in terms of the shuffle algebra on $X$.
This partly extends existing results for the lower $p$-central and $p$-Zassenhaus filtrations of $G$.
\end{abstract}

\section{Introduction}
In this paper, we describe the lower cohomology of a large family of pro-$p$ groups in terms of the combinatorics of words.
Here $p$ is a fixed prime number, and we consider the profinite cohomology groups $H^l(\bar G)=H^l(\bar G,\dbZ/p)$, $l=1,2$, where the pro-$p$ group $\bar G$ acts trivially on $\dbZ/p$.
It is well known that $H^1(\bar G)$ describes the generator structure of $\bar G$, whereas the (much deeper) group $H^2(\bar G)$ captures its relation structure \cite{NeukirchSchmidtWingberg}*{Ch.\ III, \S9}.

More specifically, let $G$ be a free pro-$p$ group on the basis $X$.
We study a category of \textsl{pro-$p$ Magnus formations}  $\Phi=(\Lam\colon G\to\dbZ_p\langle\langle X\rangle\rangle^{\times,1},\tau,e)$ -- to be described in more detail below -- which is modeled after the Magnus representation of free groups (\cite{Magnus35}, \cite{SerreGC}*{Ch.\ I, \S1.5}, \cite{Lazard54}*{Ch.\ I, \S4}).
To such a formation we associate a natural closed normal subgroup $G^\Phi$ of $G$, and take $\bar G=G/G^\Phi$.

On the combinatorial side, we consider $X$ also as an alphabet with a fixed total order.
Let $X^*$ be the set of words in $X$.
The \textsl{shuffle algebra on $X$} is the free $\dbZ$-module on $X^*$ with the \textsl{shuffle product} $\sha$ (see \S\ref{section on Lyndon words}).
Dividing by the submodule generated by all shuffle products $u\sha v$ of nonempty words $u,v$, we obtain the \textsl{indecomposable quotient} $\Sh(X)_\indec$ of the shuffle algebra.
It is this graded module which lies at the center of our description of the cohomology:
For $p$ sufficiently large, we establish a canonical isomorphism
 \begin{equation}
 \label{isomorphism}
 \bigl(\bigoplus_{s\in I}\Sh(X)_{\indec,s})\tensor(\dbZ/p)\isom H^2(G/G^\Phi),
 \end{equation}
where $I$ is a set of positive integers associated with $\Phi$ (see Theorem \ref{first isom for shuffle algebra} for the precise statement).

The isomorphism (\ref{isomorphism}) was earlier proved in some important cases:

\medskip

(1) \quad
Labute, in his seminal work \cite{Labute67} on the structure of pro-$p$ Demu\v skin groups (following Serre \cite{SerreDemuskin}), implicitly proves (\ref{isomorphism}) when $G^\Phi=G^p[G,G]$ and $I=\{1,2\}$.
Namely, the cosets of the words $(x)$, $x\in X$, and $(xy)$, where $x,y\in X$ and $x<y$, form a linear basis of the left-hand side of (\ref{isomorphism}).
Let $\varphi_x$, $x\in X$, be the basis of $H^1(G/G^\Phi)$ which is dual to $X$.
Then, for $p$ odd, the Bockstein elements $\Bock(\varphi_x)$ (see Example \ref{example with Bocksteins}) and the cup products $\varphi_x\cup\varphi_y$,  $x<y$, form a linear basis of $H^2(G/G^\Phi)$, giving the desired isomorphism.

\medskip

(2) \quad
This was extended in \cite{Efrat17} and  \cite{Efrat20} to the \textsl{lower $p$-central filtration}  $G^{(n,p)}$, $n=1,2\nek$ of $G$.
Recall that these closed subgroups of $G$ are defined inductively by
\[
G^{(1,p)}=G, \quad G^{(n+1,p)}=(G^{(n,p)})^p[G,G^{(n,p)}], \quad  n\geq1.
\]
By constructing an appropriate linear basis of $H^2(G/G^{(n,p)})$, it was shown that (\ref{isomorphism}) holds for $G^\Phi=G^{(n,p)}$, $n<p$, and $I=\{1,2\nek n\}$.
When $n=2$ this recovers Labute's result.

\medskip

(3) \quad
Let $G_{(n,p)}$, $n=1,2\nek$ be the \textsl{$p$-Zassenhaus filtration} of $G$ (also called the \textsl{modular dimension filtration} \cite{DixonDuSautoyMannSegal99}*{Ch.\ 11}; See also \cite{MinacPasiniQuadrelliTan21}, \cite{MinacPasiniQuadrelliTan22}).
Thus
\[
G_{(1,p)}=G, \quad G_{(n,p)}=(G_{(\lceil n/p\rceil,p)})^p\prod_{k+l=n}[G_{(k,p)},G_{(l,p)}], \quad n\geq2.
\]
In \cite{Efrat23} we prove (\ref{isomorphism}) when $G^\Phi=G_{(n,p)}$, $n<p$, and $I=\{1,n\}$.

\medskip

In the general situation studied in this work, the subgroups $G^\Phi$ and the isomorphism (\ref{isomorphism}) are constructed using another central notion in the combinatorics of words: \textsl{Lyndon words}.
These are the nonempty words $w$ in $X^*$ which are smaller in the alphabetical order than all their proper suffixes.
For each such word $w$, one associated its \textsl{Lie element} $\tau_w$, which is an iterated commutator in the free pro-$p$ group $G$ (Example \ref{the shuffle Magnus formation}).
Appropriate $p$-powers of the Lie elements approximately generate the subgroups $G^{(n,p)}$ and $G_{(n,p)}$ of $G$, as follows:
\begin{enumerate}
\item[(i)]
The powers $\tau_w^{p^{n-i}}$, where $w$ is a Lyndon word of length $1\leq i\leq n$, generate $G^{(n,p)}$ modulo $G^{(n+1,p)}$ \cite{Efrat17}*{Th.\ 5.3}.
\item[(ii)]
When $n\geq p$,  the powers $\tau_{(x)}^q=x^q$, where $x\in X$ and $q$ is the smallest $p$-power $\geq n$, together with the Lie elements $\tau_w$, where $w$ is a Lyndon word of length $n$, generate $G_{(n,p)}$ modulo $(G_{(n,p)})^p[G,G_{(n,p)}]$ \cite{Efrat23}*{Th.\ 4.6}.
\end{enumerate}
Such approximations seem unavailable in more general cases, and instead, we simply work with subgroups $G^\Phi$ generated by $p$-powers of the $\tau_w$.
More specifically,  we take an integer $n\geq2$ and a map $j\colon\{1,2\nek n\}\to\dbZ_{\geq0}$.
Set
\[
I=\bigl\{1\leq i\leq n\ \bigm|\ ip^{j(i)}\leq i'p^{j(i')} \hbox{ for every }1\leq i'\leq i\bigr\}.
\]
Let $L$ consist of all Lyndon words in $X^*$ with lengths in $I$, and let $\tau\colon L\to G$ be the map $w\mapsto\tau_w$.
Then the \textsl{Magnus formation}  $\Phi$ considered in (\ref{isomorphism}) is the triple $(\Lam,\tau,e)$, where $\Lam\colon G\to\dbZ_p\langle\langle X\rangle\rangle^{\times,1}$ is a continuous homomorphism into the group of $1$-units in $\dbZ_p$, and $e(i)=p^{j(i)}$.
We define $G^\Phi$ to be the closed subgroup of $G$ generated by all powers $\tau_w^{e(i)}$, $w\in L$;
See \S\ref{section on Magnus formations} for the definition of general Magnus formations.
This contains (i) and (ii) as special cases (Examples \ref{main result for lower p central filtration} and \ref{main result for p Zassenhaus filtration}).

The idea of the proof of our isomorphism theorem is to construct natural bases of the two sides of (\ref{isomorphism}):
First, using results of Radford \cite{Radford79} and, independently, Perrin and Viennot, we show that the cosets of the Lyndon words $w$ with lengths in $I$ form a linear basis of the left-hand side of (\ref{isomorphism}).
Further, to each such word $w$ we associate a cohomology element $\bar\rho_w^*(\alp_w)$ in $H^2(G/G^\Phi)$.
Using a ``triangularity" property of Lyndon words (Definition \ref{def of Magnus formations}(vi)), we show that these elements form a linear basis of $H^2(G/G^\Phi)$, which we call the \textsl{Lyndon basis}.
Associating these bases one with each other yields the desired isomorphism (\ref{isomorphism}).

The basis elements $\bar\rho_w^*(\alp_w)$ belong to an intriguing set of cohomology elements  in $H^2(G/G^\Phi)$, called the \textsl{unitriangular spectrum}.
We recall from \cite{Efrat17} that these elements are defined as follows:
Let $\dbU_i$ be the group of all unitriangular   (i.e., unipotent and upper-triangular)  $(i+1)\times(i+1)$-matrices over $\dbZ/p^{j(i)+1}$.
The natural central extension $0\to\dbZ/p\to\dbU_i\to\bar\dbU_i\to 1$ has a classifying cohomology element $\alp_i\in H^2(\bar \dbU_i)$.
For every continuous homomorphism $\bar\rho\colon G\to\bar\dbU_i$ one has the pullback $\bar\rho^*(\alp_i)$ in $H^2(G)$.
In particular, every word $w$ of length $i$ gives rise to a \textsl{Magnus  representation} $\rho_w\colon G\to\dbU_i$  (\S\ref{section on Magnus formations}), which induces a continuous homomorphism $\bar\rho_w\colon G/G^\Phi\to\bar\dbU$.
For $w$ Lyndon of length $i\in I$, we obtain the basis element $\bar\rho_w^*(\alp_i)$.

In its extreme ends, the unitriangular spectrum contains Bockstein elements (for $i=1$) and Massey product elements (when $j(i)=0$) -- see Examples \ref{example with Bocksteins} and \ref{example with Massey products}.
Bockstein elements are fairly well understood, whereas Massey product elements were extensively studied in recent years in Galois-theoretic situations (see e.g., \cite{EfratMatzri17}, \cite{GuillotMinacTopazWittenberg18}, \cite{HarpazWittenberg19}, \cite{HopkinsWickelgren15}, \cite{LamLiuSharifiWakeWang23}, \cite{MinacTan16}, as well as the references therein), and for number fields from arithmetical topology perspective (e.g., in \cite{Morishita12}, \cite{Vogel05}, \cite{KodaniMorishitaTerashima17}).
However, the behavior of the inner elements of the spectrum in such situations is still not well understood, and we hope that the connections investigated in this paper will be applicable also in these other contexts.

I thank Darij Grinberg for helpful comments on a previous version of this paper.

\section{Binomial maps}
\label{section on binomial maps}
For the rest of the paper, we fix an integer $n\geq2$.
Let $e\colon \{1,2\nek n\}\to\dbZ_{\geq1}$ be a map.

\begin{defin}
\label{def of binomiality}
\rm
The map $e$ is \textsl{binomial} if for all positive integers $i,i',l$ such that $i'l\leq i\leq n$ and $1\leq l\leq e(i')$ one has $e(i)|\binom{e(i')}l$.
\end{defin}

The following lemma shows that this definition of binomial maps coincides with the one given at \cite{ChapmanEfrat16}*{Def.\ 3.6}:

\begin{lem}
\label{equivalent conditions for binomiality}
The map $e$ is binomial if and only if the following two conditions hold:
\begin{enumerate}
\item[(a)]
For all positive integers $i',l$ such that $i'l\leq n$ and $1\leq l\leq e(i')$ one has $e(i'l)|\binom{e(i')}l$; and
\item[(b)]
For every $1\leq i'\leq i\leq n$ one has $e(i)|e(i')$.
\end{enumerate}
\end{lem}
\begin{proof}
If $e$ is binomial, then (a) and (b) hold by taking $i=i'l$ and $l=1$, respectively.
The inverse implication is immediate.
\end{proof}

We let $I_e$ be the set of all $1\leq i\leq n$ such that for every $1\leq i'\leq i$ one has
\begin{equation}
\label{def of I}
i'e(i')\geq ie(i).
\end{equation}
In particular, $1\in I_e$.

\begin{exam}
\label{first example}
\rm
Let $a_1\nek a_{n-1}$ be positive integers.
For $1\leq i\leq n$ we set
\[
e(i)=\gcd\Bigl\{\prod_{k\in K}a_k\ \Bigm|\ K\subseteq\{1,2\nek n-1\},\ |K|=n-i\Bigr\}
\]
(thus $e(n)=1$).
By \cite{ChapmanEfrat16}*{Examples 3.3, 3.8}, $e$ is binomial.

Now consider the special case where $a_k=p^{r_k}$ with $p$ a prime number and $1\leq r_1\leq\cdots \leq r_{n-1}$.
Then $e(i)=p^{j(i)}$ with $j(i)=\sum_{k=1}^{n-i}r_k$.
Since the map $x-\log_px$ is increasing in $[1,\infty)$, for $1\leq i'\leq i$ one has
\[
\log_p(i/i')\leq i-i'\leq \sum_{k=n-i+1}^{n-i'}r_k=j(i')-j(i),
\] and (\ref{def of I}) holds.
Therefore $I_e=\{1,2\nek n\}$ in this case.

In particular, when $a_k=\cdots=a_{n-1}=p$ we have $e(i)=p^{n-i}$.
\end{exam}

\begin{exam}
\label{second example}
\rm
Let $p$ be a prime number and let $t$ be a positive integer.
We set $j(i)=t\lceil\log_p(n/i)\rceil$ and $e(i)=p^{j(i)}$ for $1\leq i\leq n$.
Then $j$ is weakly decreasing, so (b) of Lemma \ref{equivalent conditions for binomiality} holds.
By \cite{ChapmanEfrat16}*{Example 3.9}, $e$ is binomial.
Also note that $j(1)\geq t$, so $e(1)>1$.

The set $I_e$ is given by the following proposition.
In the special case $t=1$ it was proved in \cite{Efrat23}*{Lemma 4.1}, and the proof in the general case is obtained by minor adjustments.
\end{exam}

\begin{prop}
\label{conditions for J(n)}
In the situation of Example \ref{second example}, $I_e$ consists of the integers $i_k=\lceil n/p^k\rceil$, where $k\geq0$.
\end{prop}
\begin{proof}
The sequence $i_k$ is weakly decreasing to $1$.
We may restrict ourselves to $k$ such that $p^k\leq n$.
Then $(n/p^k)+1\leq n/p^{k-1}$, so $n/p^k\leq i_k<n/p^{k-1}$.
Hence $j(i_k)=tk$.

Since $n/p^k\leq\lceil n/p^{k+1}\rceil p^t$, one has $i_kp^{tk}\leq i_{k+1}p^{t(k+1)}$, i.e.,  the sequence $i_kp^{j_n(i_k)}=i_kp^{tk}$ is weakly increasing in the above range.

We also observe that if $i< i_{k-1}$, then $i<n/p^{k-1}$, i.e.,  $j(i)\geq tk$.

\medskip
We now show that every $i\in I_e$ has the form $i_k$ for some $k$.
Since $n=i_0$,  we may assume that $i<n$.
Hence there is $k$ in the above range such that $i_k\leq  i< i_{k-1}$.
By the previous observation, $j(i)\geq tk$.
Taking $i'=i_k$ in (\ref{def of I}), we obtain that
\[
i_kp^{j(i_k)}\geq ip^{j(i)}\geq  i_kp^{tk}=i_kp^{j(i_k)}.
\]
Hence $i=i_k$.

Conversely, we show that each $i_k$ is in $I_e$.
Take $1\leq i'<i_k$.
There exists $l$ in the above range such that $i_l\leq i'<i_{l-1}$.
Necessarily $l>k$, so $i_lp^{tl}\geq  i_kp^{tk}$.
As we have observed, $j(i')\geq tl$.
Hence
\[
i'p^{j(i')}\geq  i_lp^{tl}\geq  i_kp^{tk}=i_kp^{j(i_k)}.
\qedhere
\]
\end{proof}

\section{Words}
\label{section on Lyndon words}
We refer to \cite{Reutenauer93}, \cite{Lothaire83}, and \cite{GrinbergReiner20} as general sources on the combinatorics of words.

We fix a nonempty totally ordered  set $(X,\leq)$, considered as an alphabet.
Let $X^*$ be the free unital monoid on $X$ with the concatenation product.
Its elements are considered as (associative) words in $X$, and the unit element $1$ of $X^*$ is the empty word.
We write $|w|$ for the length of the word $w$.
Let $\leq_{\textrm{alp}}$ be the alphabetical (lexicographic) total order on $X^*$.
The \textsl{length-alphabetical} total order $\preceq$ on $X^*$ is defined by $w_1\preceq w_2$ if and only if $|w_1|<|w_2|$, or both $|w_1|=|w_2|$ and $w_1\leq_{\textrm{alp}}w_2$.

A nonempty word $w\in X^*$ is called a \textsl{Lyndon word} if it is smaller in $\leq_{\textrm{alp}}$ than all its proper right factors (i.e., suffixes).
See \cite{Deligne10}*{\S2.1} for an equivalent definition.

Every Lyndon word $w$ of length $\geq2$ has a \textsl{standard factorization} as a concatenation $w=uv$ of Lyndon words $u,v$.
Namely, $v$ is the $\leq_{\textrm{alp}}$-minimal nontrivial right factor of $w$ which is a Lyndon word;
equivalently, $v$ is the longest nontrivial right factor of $w$ which is a Lyndon word \cite{Efrat23}*{Lemma 2.2}.
The set of all Lyndon words in $X^*$ is a Hall set \cite{Reutenauer93}*{Th.\ 5.1}, and the above factorization coincides with the general notion of a standard factorization in Hall sets \cite{Reutenauer93}*{\S4.1}.

Consider a unital commutative ring $R$.
Let $R\langle\langle X\rangle\rangle$ denote the $R$-module of formal power series $f$ with coefficients in $R$ and over the set $X$ of non-commuting variables.
Thus $f=\sum_{w\in X^*}f_ww$ with $f_w\in R$.
The concatenation induces on $R\langle\langle X\rangle\rangle$ via linearity the structure of an $R$-algebra.
Let $R\langle X\rangle$ be the subalgebra of $R\langle\langle X\rangle\rangle$ consisting of all noncommutative polynomials in $X$, that is, all such power series $f$ whose support
\[
\Supp(f)=\{w\in X^*\ |\ f_w\neq0\}
\]
is finite.
We identify a word $w\in X^*$ as a monomial in $R\langle X\rangle$.

For $f=\sum_{w\in X^*}f_ww\in R\langle\langle X\rangle\rangle$ and $g=\sum_{w\in X^*}g_ww\in R\langle X\rangle$ we define the \textsl{scalar product}
\begin{equation}
\label{scalar product}
(f,g)_R=\sum_{w\in X^*}f_wg_w.
\end{equation}
It is $R$-bilinear, and $(f,w)_R=f_w$ and $(w,g)_R=g_w$ for every $w$.

Next let $\star\colon X^*\times X^*\to R\langle X\rangle$ be a binary map such that
the module $R\langle X\rangle$ is a unital associative $R$-algebra with respect to the induced $R$-bilinear map
\[
\star \colon R\langle X\rangle\times R\langle X\rangle\to R\langle X\rangle, \ \
\bigl(\sum_uf_uu\bigr)\star\bigl(\sum_vg_vv\bigr)=\sum_{u,v} f_ug_v(u\star v),
\]
and the unit element $1\in X^*$.
We further assume that for nonempty $u,v\in X^*$ the support $\Supp(u\star v)$ consists only of words $w$ of length  $1\leq|w|\leq|u|+|v|$.
We say that $f\in R\langle\langle X\rangle\rangle$ is \textsl{compatible with $\star$} if for every nonempty words $u,v\in X^*$ one has
\begin{equation}
\label{generalized CFL}
(f,u)_R\cdot (f,v)_R=(f,u\star v)_R.
\end{equation}

The proof of the following fact is straightforward.
Here $(u\star v)_s$ denotes the homogeneous part of degree $s$ of $u\star v$.

\begin{lem}
\label{series satisfying CFL relation}
Suppose that $f\in R\langle\langle X\rangle\rangle$ is compatible with $\star$, and let $N$ be an ideal in $R$.
Let $u,v$ be nonempty words in $X^*$ with $s=|u|+|v|$.
If $(f,w)_R\in N$ for every nonempty $w\in X^*$ with $|w|<s$, then $(f,(u\star v)_s)_R\in N$.
\end{lem}

We will be especially interested in the \textsl{shuffle product} $u\sha v$ and \textsl{infiltration product} $u\downarrow v$ of words $u,v\in X^*$  defined as follows (see \cite{ChenFoxLyndon58},  \cite{Reutenauer93}*{pp.\ 134--135}):
Write $u=(x_1\cdots x_r)$, $v=(x_{r+1}\cdots x_{r+t})\in  X^*$.
Then
\[
u\sha v=\sum_\sig(x_{\sig\inv(1)}\cdots x_{\sig\inv(r+t)})\in R\langle X\rangle,
\]
where $\sig$ ranges over all permutations of $1,2\nek r+t$ such that $\sig(1)<\cdots<\sig(r)$ and $\sig(r+1)<\cdots<\sig(r+t)$.

Similarly, consider all surjective maps $\sig\colon\{1,2\nek r+t\}\to\{1,2\nek k(\sig)\}$, with $1\leq k(\sig)\leq r+t$, $\sig(1)<\cdots<\sig(r)$ and $\sig(r+1)<\cdots<\sig(r+t)$, and which satisfy the following
weak form of injectivity:
If $\sig(i)=\sig(j)$, then $x_i=x_j$.
Then we set
\begin{equation}
\label{infiltration}
u\downarrow v=\sum_\sig(x_{\sig\inv(1)}\cdots x_{\sig\inv(k(\sig))})\in R\langle X\rangle.
\end{equation}
By our assumption, $x_{\sig\inv(i)}$ does not depend on the choice of the preimage $\sig\inv(i)$.
Thus  $u\sha v$ is the part of $u\downarrow v$ of degree $r+t$, that is, the partial sum corresponding to all such maps $\sig$ which are bijective.

The maps $R\langle X\rangle\times R\langle X\rangle\to R\langle X\rangle$ induced by  $\sha$ and $\downarrow$ as above are commutative and associative
\cite{Lothaire83}*{p.\ 128, Prop.\ 6.3.15}, and the assumption about $\Supp(u\star v)$ is clearly satisfied.

\section{Magnus formations}
\label{section on Magnus formations}
The discussion in this section can be carried out both in the context of discrete groups and rings, as well as profinite groups and rings.
Since we are motivated by profinite applications, and the discrete setting is easily obtained from the profinite one by obvious amendments, we carry the discussion in the latter setting.

Let $R$ be a profinite commutative ring with unit $1_R$.
We equip $R\langle\langle X\rangle\rangle$ with the minimal topology for which the maps $(\cdot, w)_R\colon R\langle\langle X\rangle\rangle\to R$ are continuous for every $w\in X^*$.
Then $R\langle\langle X\rangle\rangle$ is also a profinite ring, and both its subgroup $R\langle\langle X\rangle\rangle^\times$ of invertible elements and its subgroup $R\langle\langle X\rangle\rangle^{\times,1}$ of invertible elements $f$ with $(f,1)_R=1_R$ are profinite groups \cite{Efrat14}*{\S5}.

Let $G$ be a profinite group, and consider a continuous map
\[
\Lam\colon G\to R\langle\langle X\rangle\rangle^{\times,1}, \quad
\Lam(\sig)=\sum_{w\in X^*}\eps_w(\sig)w.
\]
The maps $\eps_w\colon G\to R$, $\eps(\sig)=(\Lam(\sig),w)_R$, are then continuous.

Let $\dbU_i(R)$ be the group of all $(i+1)\times(i+1)$-matrices over $R$ which are \textsl{unitriangular}, i.e., unipotent and upper-triangular.
It is a profinite group with respect to the natural (product) topology induced from $R$.
For a word $w=(a_1a_2\cdots a_i)\in X^*$ of length $i$ we define a map
\[
\rho_{w,R}\colon G\to \dbU_i(R), \quad
\rho_{w,R}(\sig)_{kl}=\eps_{(a_ka_{k+1}\cdots a_{l-1})}(\sig)
\]
for $1\leq k\leq l\leq i+1$.
It is continuous, and we call it the \textsl{Magnus representation} of $G$ over $R$ associated with $w$.

\begin{prop}
\label{Lambda a homomorphism}
The following conditions are equivalent:
\begin{enumerate}
\item[(a)]
$\Lam$ is a group homomorphism;
\item[(b)]
For every $w\in X^*$ and $\sig_1,\sig_2\in G$ one has
\[
\eps_w(\sig_1\sig_2)=\sum_{u,v\in X^*,\ w=uv}\eps_u(\sig_1)\eps_v(\sig_2),
\]
where the summation is over all decompositions of $w$ as a concatenation $w=uv$.
\item[(c)]
For every $w\in X^*$ the map $\rho_{w,R}$ is a group homomorphism.
\end{enumerate}
\end{prop}
\begin{proof}
The implications (a)$\Leftrightarrow$(b)$\Rightarrow$(c) are straightforward.
For (c)$\Rightarrow$(b) we look at the $(1,i+1)$-entry of $\rho_{w,R}$.
\end{proof}

Given $w\in X^*$, we write $o(w)$ for a general element of $R\langle\langle X\rangle\rangle$ whose support consists of words strictly larger than $w$ with respect to the length-alphabetical order $\preceq$ (see \S\ref{section on Lyndon words}), and composed of letters appearing in $w$.

The  \textsl{lower central series} $G^{(i)}$, $i=1,2\nek$ of $G$ is defined inductively by $G^{(1)}=G$ and $G^{(i+1)}=[G,G^{(i)}]$ (in the profinite sense).
We recall that $n\geq2$ is a fixed integer as in \S\ref{section on binomial maps}.

\begin{defin}
\label{def of Magnus formations}
\rm
A (profinite) \textsl{Magnus formation} over $R$ is a triple
\[
\Phi=\bigl(\Lam\colon G\to R\langle\langle X\rangle\rangle^{\times,1},\tau\colon L\to G,e\colon \{1,2\nek n\}\to\dbZ_{\geq1}\bigr),
\]
such that:
\begin{enumerate}
\item[(i)]
$G$ is a profinite group;
\item[(ii)]
$\Lam$ is a continuous group homomorphism;
\item[(iii)]
$L$ is a nonempty subset of $X^*$;
\item[(iv)]
The words in $L$ have lengths in $\{1,2\nek n\}$;
\item[(v)]
$\tau$ is a map such that $\tau(w)\in G^{(i)}$ for every $w\in L$ of length $i$;
\item[(vi)]
For every $w\in L$ one has $\Lam(\tau(w))=1+w+o(w)$ (\textsl{triangularity});
\item[(vii)]
The map $e$ is binomial and is not identically $1$.
\end{enumerate}
\end{defin}

Consider Magnus formations
\[
\Phi_l=(\Lam_l\colon G_l\to R\langle\langle X\rangle\rangle^{\times,1},\tau_l\colon L_l\to G_l,e), \quad l=1,2,
\]
over $R$, where $L_1\subseteq L_2$.
A \textsl{morphism} $\Phi_1\to\Phi_2$  is a continuous group homomorphism $\gam\colon G_1\to G_2$ such that $\Lam_1=\Lam_2\circ\gam$ and $\tau_2|_{L_1}=\gam\circ\tau_1$.

Let $\star\colon X^*\times X^*\to R\langle X\rangle$ be a binary map as in \S\ref{section on Lyndon words}.
We say that the Magnus formation $\Phi=(\Lam,\tau,e)$ is \textsl{compatible with $\star$} if for every $\sig\in G$ the series $\Lam(\sig)$ is compatible with $\star$, in the sense of \S\ref{section on Lyndon words}.

We note that if $\gam\colon\Phi_1\to\Phi_2$ is a morphism of Magnus formations as above, and if $\Phi_2$ is compatible with $\star$, then $\Phi_1$ is also compatible with $\star$.
The converse holds if $\gam\colon G_1\to G_2$ is surjective.

\begin{exam}
\label{the shuffle Magnus formation}
\rm
Let $p$ be a prime number, let $G$ be the free pro-$p$ group on the alphabet $X$ as a basis \cite{FriedJarden08}*{\S17.4}, and take $R=\dbZ_p$.
The \textsl{pro-$p$ Magnus homomorphism} $\Lam_{\dbZ_p}\colon G\to \dbZ_p\langle\langle X\rangle\rangle^{\times,1}$ is defined on the free generators $x\in X$ of $G$ by $\Lam_{\dbZ_p}(x)=1+x$ (note that $(1+x)\sum_{k\geq0}(-1)^kx^k=1$, so indeed $1+x\in  \dbZ_p\langle\langle X\rangle\rangle^{\times,1}$).

For a Lyndon word $w$ in $X^*$ one defines its \textsl{Lie element} $\tau_w$ in $G$ by induction on $|w|$ as follows:
When $w=(x)$ has length $1$, set $\tau_{(x)}=x$.
Otherwise let $w=uv$ be the standard factorization of $w$ (see \S\ref{section on Lyndon words}), and set $\tau_w=[\tau_u,\tau_v]$.

Let $e\colon \{1,2\nek n\}\to\dbZ_{\geq1}$ be any binomial map which is not identically $1$, and let $L$ be a nonempty set of Lyndon words in $X^*$ of lengths $\leq n$.
We set $\tau\colon L\to S$ to be the map $w\mapsto \tau_w$.
Then (v) holds by \cite{Reutenauer93}*{Cor.\ 6.16}, and (vi) holds by \cite{Efrat17}*{Prop.\ 4.4(c)}.
Consequently, $(\Lam_{\dbZ_p},\tau,e)$ is a Magnus formation over $\dbZ_p$.

Moreover, this formation is compatible with the infiltration product $\downarrow$:
Indeed, this was shown by Chen, Fox, and Lyndon in the discrete case for $R=\dbZ$ \cite{ChenFoxLyndon58}*{Th.\ 3.6}, and the case $R=\dbZ_p$ follows (see also \cite{Morishita12}*{Prop.\ 8.6}, \cite{Reutenauer93}*{Proof of Th.\ 6.4}, \cite{Vogel05}*{Prop.\ 2.25}).
\end{exam}

\section{The fundamental matrix}
We consider a Magnus formation
\[
\Phi=(\Lam\colon G\to R\langle\langle X\rangle\rangle^{\times,1},\tau\colon L\to G,e\colon \{1,2\nek n \}\to\dbZ_{\geq1})
\]
over the profinite ring $R$.
For $w\in L$ of length $i$ let
\[
\sig_w=\tau(w)^{e(i)}.
\]

\begin{prop}
\label{p adic value of Magnus coefficient}
Let  $w\in X^*$ and $w'\in L$ have lengths $i,i'$, respectively.
\begin{enumerate}
\item[(a)]
If $1\neq w\prec w'$, then $\eps_w(\sig_{w'})=0$;
\item[(b)]
If $w=w'$, then $\eps_w(\sig_{w'})=e(i)1_R$;
\item[(c)]
If $1\leq i\leq n$, then $\eps_w(\sig_{w'})\in e(i)R$;
\item[(d)]
If $w$ has letters not appearing in $w'$, then $\eps_w(\sig_{w'})=0$.
\end{enumerate}
\end{prop}
\begin{proof}
By the triangularity property (Definition \ref{def of Magnus formations}(vi)) and the binomial expansion formula,
\[
\Lam(\sig_{w'})=\Lam(\tau(w'))^{e(i')}=(1+w'+o(w'))^{e(i')}=1+e(i')w'+o(w').
\]
This gives (a) and (b).
Moreover,
\[
\Lam(\sig_{w'})=\Lam(\tau(w'))^{e(i')}\in\biggl(1+\sum_{|u|\geq i'}Ru\biggr)^{e(i')}\subseteq
1+\sum_{1\leq l\leq e(i'), i'l\leq|v|}\binom{e(i')}l Rv.
\]
If $l$ satisfies $1\leq l\leq e(i')$ and $i'l\leq i$, then the binomiality of $e$ implies that $e(i)|\binom{e(i')}l$.
Taking in the right sum $v=w$, we deduce (c).

Moreover, in the sums above we may restrict to words $u$, $v$ whose letters appear in $w'$.
This shows (d).
\end{proof}

Now fix an integer $m\geq2$ such that for every positive integer $e$ one has $\dbZ/m\dbZ\isom eR/meR$ via the group homomorphism $k\mapsto ke1_R$.
In the applications, we will take $m=p$ prime and $R=\dbZ_p$,  and this condition clearly holds.

For $1\leq i\leq n$ let $\pi_i\colon R\to R/me(i)R$ be the natural epimorphism.
For $w\in X^*$ and $w'\in L$ of lengths $1\leq i,i'\leq n$, respectively, we set
\[
\langle w,w'\rangle=\pi_i(\eps_w(\sig_{w'})).
\]
By Proposition \ref{p adic value of Magnus coefficient}(c), $\langle w,w'\rangle\in e(i)R/me(i)R$.
Under the identification $\dbZ/m\isom e(i)R/me(i)R$, we may consider $\langle w,w'\rangle$ as an element of $\dbZ/m$.

We call the transposed (possibly infinite) matrix
\[
\Bigl[\langle w,w'\rangle\Bigr]_{w,w'\in L}^T,
\]
where $L$ is totally ordered by $\preceq$, the \textsl{fundamental matrix} of the Magnus formation $\Phi$.
From Proposition \ref{p adic value of Magnus coefficient}(a)(b) we deduce:

\begin{cor}
\label{fundamental matrix is unitriangular}
The fundamental matrix of a Magnus formation is unitriangular.
\end{cor}

\begin{rem}
\label{the fundamental matrix is functorial}
\rm
This construction is functorial in the following sense:
Consider Magnus formations
\[
\Phi_l=(\Lam_l\colon G_l\to R\langle\langle X\rangle\rangle^{\times,1},\tau_l\colon L_l\to G_l,e\colon \{1,2\nek n\}\to\dbZ_{\geq1}), \  l=1,2,
\]
over $R$, with $L_1\subseteq L_2$.
Let $\gam\colon \Phi_1\to\Phi_2$ be a morphism as in \S\ref{section on Magnus formations}.
For every $w'\in L_1$ of length $i'$ we have
\[
\Lam_1(\tau_1(w')^{e(i')})=\Lam_2(\tau_2(w')^{e(i')}).
\]
Hence the values of $\langle w,w'\rangle$ with respect to $\Phi_1$ and to $\Phi_2$ coincide.
\end{rem}

\section{Unitriangular representations}
\label{section on unitriangular representations}
For an integer $i\geq1$ and a profinite commutative unital ring $R$, we write $\Id$ for the identity matrix in the group $\dbU_i(R)$ of unitriangular $(i+1)\times(i+1)$-matrices over $R$ (see \S\ref{section on Magnus formations}).
We also write $E_{kl}$ for the $(i+1)\times(i+1)$-matrix over $R$ which is $1$ at entry $(k,l)$, and is $0$ elsewhere.

From now on we assume that $m=p$ is a prime number and $R=\dbZ_p$.

We record the following fact about the lower central series of $\dbU_i$, proved in  \cite{Efrat23}*{Prop.\ 6.2(c)}:

\begin{lem}
\label{lemma from Jussieu paper}
For integers $i\geq i'\geq1$ and $j,j'\geq0$, one has  $j'\geq j+\log_p(i/i')$ if and only if $(\dbU_i(\dbZ/p^{j+1})^{(i')})^{p^{j'}}\leq \Id+p^j\dbZ E_{1,i+1}$.
\end{lem}

Let $\Phi=(\Lam\colon G\to \dbZ_p\langle\langle X\rangle\rangle^{\times,1},\tau\colon L\to G,e)$ be a Magnus formation over $\dbZ_p$, where $e(i)=p^{j(i)}$ for some map $j\colon \{1,2\nek n\}\to\dbZ_{\geq0}$.
Let $I_e$ be the subset of $\{1,2\nek n\}$ defined in \S\ref{section on binomial maps}.

Given a word $w\in X^*$ of length $1\leq i\leq n$ we set
\[
R_w=\dbZ/p^{j(i)+1}(=\dbZ/me(i)), \qquad
\dbU_w=\dbU_i(R_w).
\]
Let $\dbU_w^0$ be the subgroup $\Id+p^{j(i)}\dbZ E_{1,i+1}$ of $\dbU_w$.
Then $\dbU_w^0\isom\dbZ/p$, and $\dbU_w^0$ is central in $\dbU_w$.
One has $\dbU_w^0=(\dbU_w^{(i)})^{p^{j(i)}}$ \cite{Efrat23}*{Prop.\ 6.2(b)}.

Let $\bar\Lam_w\colon G\to R_w\langle\langle X\rangle\rangle^{\times,1}$ be the continuous homomorphism induced by $\Lam$, and let $\rho_{w,R_w}\colon G\to\dbU_w$ be the continuous Magnus representation associated with $\bar\Lam_w$,  as in \S\ref{section on Magnus formations}.

We define $G^\Phi$ to be the closed normal subgroup of $G$ generated by all powers $\sig_w=\tau(w)^{p^{j(i)}}$, where $w\in L$ and $i=|w|$.

\begin{prop}
\label{image of SL}
\begin{enumerate}
\item[(a)]
For a word $w\in X^*$ of length $i\in I_e$ one has $\rho_{w,R_w}(G^\Phi)\leq\dbU^0_w$.
\item[(b)]
If in addition $w\in L$, then $\rho_{w,R_w}(G^\Phi)=\dbU^0_w$.
\end{enumerate}
\end{prop}
\begin{proof}
(a) \quad
Since $\dbU^0_w$ is normal in $\dbU_w$, it suffices to show that $\rho_{w,R_w}(\sig_{w'})\in\dbU^0_w$ for every $w'\in L$.
Let $i'=|w'|$.

Assume first that $i'\leq i$.
Since $i\in I_e$, (\ref{def of I}) holds, so  $j(i')\geq j(i)+\log_p(i/i')$.
Hence, by Lemma \ref{lemma from Jussieu paper}, $(\dbU_w^{(i')})^{p^{j(i')}}\leq  \Id+p^{j(i)}\dbZ E_{1,i+1}=\dbU^0_w$.
By condition (v) of Definition \ref{def of Magnus formations},  $\tau(w')\in G^{(i')}$, and we deduce that
\[
\rho_{w,R_w}(\sig_{w'})=\rho_{w,R_w}(\tau(w'))^{p^{j(i')}}\in(\dbU_w^{(i')})^{p^{j(i')}}\leq\dbU^0_w.
\]

When $i<i'$ Proposition \ref{p adic value of Magnus coefficient}(a) implies that $\rho_{w,R_w}(\sig_{w'})=\Id$.

\medskip

(b) \quad
By Proposition \ref{p adic value of Magnus coefficient}(b),
$\rho_{w,R_w}(\sig_w)_{1,i+1}=p^{j(i)}1_{R_w}$.
In view of (a), $\rho_{w,R_w}(\sig_w)=\Id+p^{j(i)}E_{1,i+1}$, which is a generator of $\dbU^0_w$.
\end{proof}

By Proposition \ref{image of SL}(a), for $w$ of length $i\in I_e$, the representation $\rho_{w,R_w}$ induces a continuous homomorphism
\[
\rho_w^0\colon G^\Phi\to \dbU^0_w.
\]
For every $w,w'\in L$ with $i=|w|\in I_e$ we have, under the identifications  $\dbU_w^0=p^{j(i)}\dbZ/p^{j(i)+1}\dbZ=\dbZ/p$, that
\begin{equation}
\label{pairing ww and rho sig}
\rho_w^0(\sig_{w'})=\pi_i(\eps_w(\sig_{w'}))=\langle w,w'\rangle.
\end{equation}
Since $\dbU_w^0$ is central in $\dbU_w$, the homomorphism $\rho_w^0$ is $G$-invariant, i.e., $\rho_w^0(\lam\sig\lam \inv)=\rho_w^0(\sig)$ for every $\sig\in G^\Phi$ and $\lam\in G$.

Let $\Hom(G^\Phi,\dbZ/p)$ denote the group of all continuous homomorphisms $\psi\colon G^\Phi\to \dbZ/p$.
Let $\Hom(G^\Phi,\dbZ/p)^G$ be its subgroup consisting of all such homomorphisms which are $G$-invariant.

We will need the following fact in linear algebra \cite{Efrat17}*{Lemma 8.4}:

\begin{lem}
\label{linear algebra lemma}
Let $R$ be a commutative ring and let $(\cdot,\cdot)\colon A\times B\to R$ be a non-degenerate bilinear map of $R$-modules.
Let $(\calL,\leq)$ be a finite totally ordered set, and for every $w\in \calL$ let $a_w\in A$, $b_w\in B$.
Suppose that the matrix $\bigl((a_w,b_{w'})\bigr)_{w,w'\in \calL}$ is invertible, and that $a_w$, $w\in \calL$, generate $A$.
Then $a_w$, $w\in \calL$, is an $R$-linear basis of $A$, and $b_w$, $w\in \calL$, is an $R$-linear basis of $B$.
\end{lem}

Specifically, for $\Phi$ as above, there is a well defined perfect bilinear map of $\dbZ/p$-linear spaces
\begin{equation}
\label{perfect pairing}
(\cdot,\cdot)\colon
G^\Phi/(G^\Phi)^p[G,G^\Phi]\times\Hom(G^\Phi,\dbZ/p)^G\to\dbZ/p,\quad (\bar\sig,\psi)=\psi(\sig)
\end{equation}
\cite{EfratMinac11}*{Cor.\ 2.2}.

\begin{prop}
\label{base}
Suppose that the lengths of the words in $L$ are in $I_e$.
Then the maps $\rho_w^0$, $w\in L$, form a $\dbZ/p$-linear basis of $\Hom(G^\Phi,\dbZ/p)^G$.
\end{prop}
\begin{proof}
Assume first that $L$ is finite.
The left direct factor in (\ref{perfect pairing}) is generated by the cosets of $\sig_{w'}$, $w'\in L$.
By (\ref{pairing ww and rho sig}), the matrix $[(\bar\sig_{w'},\rho^0_w)]_{w,w'\in L}$ is the transpose of the fundamental matrix of $\Phi$.
By Corollary \ref{fundamental matrix is unitriangular}, it is invertible.
Lemma \ref{linear algebra lemma} therefore implies the assertion in this case.

In the general case, write $L=\bigcup_\alp L_\alp$, where the $L_\alp$ form a direct system of finite subsets of $L$.
For every $\alp$ let $\Phi_\alp=(\Lam,\tau|_{L_\alp},e)$ be the restricted Magnus formation.
Thus $G^{\Phi_\alp}$ be the closed normal subgroup  of $G$ generated by the $\sig_w$, $w\in L_\alp$.
Then $G^\Phi$ is generated by the $G^{\Phi_\alp}$, and it follows that $\Hom(G^\Phi,\dbZ/p)^G=\dirlim_\alp\Hom(G^{\Phi_\alp},\dbZ/p)^G$.
The assertion therefore follows from the finite case.
\end{proof}

\section{The isomorphism theorem}
As before, let $m=p$ be a prime number, let $R=\dbZ_p$,  and let $\Phi=(\Lam\colon G\to\dbZ_p\langle\langle X\rangle\rangle^{\times,1},\tau\colon L\to G,e)$ be a Magnus formation over $\dbZ_p$.
We now assume further that $\Phi$ is compatible with a binary operation $\star$ as in \S\ref{section on Magnus formations}.
We denote the set of all words of length $s$ in $X^*$ by $X^s$.
For $f\in \dbZ_p\langle X\rangle$ let $f_s$ be its homogenous part of degree $s$.
Let $I_e$ be as in \S\ref{section on binomial maps}.
Recall that for $w\in X^*$ of length in $I_e$, we may view $\rho_w^0$ as an element of $\Hom(G^\Phi,\dbZ/p)^G$.

\begin{lem}
\label{generic shuffle relation}
For every nonempty words $u,v\in X^*$ with $s=|u|+|v|\in I_e$, one has
\[
\sum_{w\in X^s}(u\star v,w)_{\dbZ_p}\rho^0_w=0.
\]
\end{lem}
\begin{proof}
First, we note that since $u\star v\in\dbZ_p\langle X\rangle$, the sum is well defined.

Let $w'\in L$.
For $w\in X^*$ of length  $1\leq i<s$ we have $ip^{j(i)}\geq sp^{j(s)}$, by (\ref{def of I}), whence $j(i)>j(s)$.
By Proposition \ref{p adic value of Magnus coefficient}(c),
$\eps_w(\sig_{w'})\in p^{j(i)}\dbZ_p\subseteq p^{j(s)+1}\dbZ_p$.
We apply Lemma \ref{series satisfying CFL relation} for the ideal $p^{j(s)+1}\dbZ_p$ of $\dbZ_p$ and $f=\Lam(\sig_{w'})$ to deduce that
\[
(\Lam(\sig_{w'}),(u\star v)_s)_{\dbZ_p}
\in p^{j(s)+1}\dbZ_p.
\]
Therefore for the substitution pairing  (\ref{perfect pairing}) we have, using (\ref{pairing ww and rho sig}),
\[
\begin{split}
&(\bar\sig_{w'},\sum_{w\in X^s}(u\star v,w)_{\dbZ_p}\rho^0_w)
=\sum_{w\in X^s}(u\star v,w)_{\dbZ_p}\rho_w^0(\bar\sig_{w'})\\
&=\sum_{w\in X^s}(u\star v,w)_{\dbZ_p}\pi_s(\eps_w(\sig_{w'}))
=\pi_s\biggl(\sum_{w\in X^s}(u\star v,w)_{\dbZ_p}\eps_w(\sig_{w'})\biggr)\\
&=\pi_s\bigl((\Lam(\sig_{w'}),(u\star v)_s)_{\dbZ_p}\bigr)=0.
\end{split}
\]
Since the bilinear map (\ref{perfect pairing}) is non-degenerate, this gives the assertion.
\end{proof}

Suppose that $A=\bigoplus_{s\geq0}A_s$ is a graded $R$-module, which is a (not necessarily graded) associative $R$-algebra with respect to a product map $\circ$.
Let $N$ be the $R$-submodule of $A$ generated by the homogenous parts $(a\circ b)_{r+t}$ of the products $a\circ b$, where $a,b$ are homogenous elements of $A$ of degrees $r,t\geq1$, respectively.
It is a graded submodule, giving rise to a graded quotient $R$-module $A_\indec=A/N$, called the \textsl{indecomposable quotient} of $A$.
Let  $A_{\indec,s}$ be the homogenous component of $A_\indec$ of degree $s$.

Now take the $\dbZ_p$-module $A=\bigoplus_{w\in X^*}\dbZ_pw$ with the product map $\star$ as in \S\ref{section on Lyndon words}.
We write $A_\indec^{(L)}$ for the submodule of $A_\indec$ generated by the images $\bar w$ of the words $w$ in $L$.

\begin{thm}
\label{indecomposable quotient as the space of maps}
Let $\Phi$ be a Magnus formation over $\dbZ_p$ as above which is compatible with $\star$.
Suppose that the words in $L$ have lengths in $I_e$.
Then the map $\bar w\tensor1\mapsto\rho^0_w$,  $w\in L$, induces an isomorphism of $\dbZ/p$-linear spaces
\[
A_{\indec}^{(L)}\tensor(\dbZ/p)\xrightarrow{\ \sim\ }  \Hom(G^\Phi,\dbZ/p)^G.
\]
\end{thm}
\begin{proof}
There is a unique  $\dbZ_p$-module homomorphism
\[
h\colon\bigoplus_{s\in I_e}\bigoplus_{w\in X^s}\dbZ_pw\to \Hom(G^\Phi,\dbZ/p)^G.
\]
such that $h(w)=\rho_w^0$ for $s\in I_e$ and $w\in X^s$.
By Lemma \ref{generic shuffle relation}, $h$ is trivial on the homogenous components $(u\star v)_s$, where $u,v$ are nonempty words and $s=|u|+|v|\in I_e$.
Therefore $h$ factors via $\bigoplus_{s\in I_e}A_{\indec,s}$.
Since the lengths of the words in $L$ are in $I_e$, the homomorphism $h$ induces a $\dbZ/p$-linear map
\[
\bar h\colon A^{(L)}_\indec\tensor(\dbZ/p)\to \Hom(G^\Phi,\dbZ/p)^G,
\]
where $\bar h(\bar w\tensor1)=\rho^0_w$ for $w\in L$.

Now the $\bar w\tensor1$, where $w\in L$, span $A_{\indec}^{(L)}\tensor(\dbZ/p)$ as a $\dbZ/p$-linear space.
Furthermore, by Proposition \ref{base}, their images $\rho^0_w$ form a linear basis of the $\dbZ/p$-linear space $ \Hom(G^\Phi,\dbZ/p)^G$.
Therefore $\bar h$ is an isomorphism.
\end{proof}

\begin{rem}
\rm
This isomorphism is functorial in the following sense:
Consider the setup of Remark \ref{the fundamental matrix is functorial} (with $m=p$ and $R=\dbZ_p$), and let $A=\bigoplus_{w\in X^*}\dbZ_pw$ be as above.
Suppose that the Magnus formation $\Phi_2$ (whence also the formation $\Phi_1$) is compatible with the product map $\star$.
Then $\gam$ induces a commutative square
\[
\xymatrix{
A_{\indec}^{(L_1)}\tensor(\dbZ/p)\ar[d]\ar[r]^{\sim\ \ }& \Hom(G^{\Phi_1},\dbZ/p)^G\\
A_{\indec}^{(L_2)}\tensor(\dbZ/p)\ar[r]^{\sim\ \ }& \Hom(G^{\Phi_2},\dbZ/p)^G,\ar[u]
}
\]
where the right vertical map is the restriction.
\end{rem}

\section{The shuffle algebra}
\label{section on shuffle algebra}
Our main examples concern the $p$-adic Magnus formations of Example \ref{the shuffle Magnus formation}, in connection with the shuffle algebra.

Every word $w\in X^*$ can be uniquely written as a concatenation $w=u_1^{k_1}\cdots u_t^{k_t}$, where $u_1>_{\mathrm{alp}}\cdots>_{\mathrm{alp}}u_t$ are Lyndon words in $X^*$, $k_1>\cdots>k_t$, and where $u^k$ denotes the concatenation of $u$ with itself $k$ times \cite{Reutenauer93}*{ Cor.\ 4.7 and Th.\ 5.1}.
Consider the noncommutative polynomials
\[
Q_w=\dfrac1{k_1!\cdots k_t!}u_1^{\sha k_1}\sha\cdots\sha u_t^{\sha k_t}\in\dbQ\langle X\rangle,
\]
where $u^{\sha k}$ denotes the $k$-times shuffle product $u\sha\cdots\sha u$.
The polynomial $Q_w$ is homogenous of degree $|w|=k_1|u_1|+\cdots+k_t|u_t|$, and in fact, $Q_w\in\dbZ\langle X\rangle$ \cite{Reutenauer93}*{Th.\ 6.1}.
By a result of Radford \cite{Radford79} and Perrin and Viennot (unpublished) -- see again \cite{Reutenauer93}*{Th.\ 6.1} --
for every word $w\in X^*$ one has
\[
Q_w=w+\sum_{v\leq_{\textrm{alp}}w}a_{v,w}v
\]
for some nonnegative integers $a_{v,w}$, and where for all but finitely many $v$ we have $a_{v,w}=0$.
We may restrict to $v$ such that $|v|=|w|$.
Hence for every positive integer $s$ we have
\begin{equation}
\label{Radford basis}
\bigoplus_{w\in X^s}\dbZ Q_w=\bigoplus_{w\in X^s}\dbZ w.
\end{equation}

We recall that the \textsl{shuffle $\dbZ$-algebra} $\Sh(X)$ over $X$ is the $\dbZ$-module $\bigoplus_{w\in X^*}\dbZ w$ with the shuffle product $\sha$ (\S\ref{section on Magnus formations}).
Let $\Sh(X)_\indec$ be its indecomposable quotient with respect to the product map $\circ=\sha$.
Since for words $u,v$ with $s=|u|+|v|$ we have $(u\downarrow v)_s=u\sha v$, the indecomposable quotient of $\bigoplus_{w\in X^*}\dbZ w$ with respect to $\circ=\downarrow$ is also $\Sh(X)_\indec$.
The shuffle product $\sha$ extends in an obvious way to the $\dbZ_p$-module $\bigoplus_{w\in X^*}\dbZ_pw$, giving rise to a $\dbZ_p$-shuffle algebra $\Sh(X)\tensor\dbZ_p$.
We note that
\begin{equation}
\label{indec commutes with tensor product}
\Sh(X)_\indec\tensor\dbZ_p=(\Sh(X)\tensor\dbZ_p)_\indec.
\end{equation}

\begin{prop}
\label{the tensorized shuffle module}
Let $L$ be the set of all Lyndon words in $X^*$ of length $s\geq1$.
Let $r$ be a positive integer whose prime factors are larger than $s$.
Then
\[
\Sh(X)_{\indec,s}\tensor(\dbZ/r)=\Sh(X)_{\indec,s}^{(L)}\tensor(\dbZ/r).
\]
\end{prop}
\begin{proof}
Take $w\in X^s$ with its decomposition $w=u_1^{k_1}\cdots u_t^{k_t}$ as above.
As $k_1\nek k_t\leq s$, the assumption on $r$ implies that $k_1!\cdots k_t!$ is invertible in $\dbZ/r$.
If in addition $w\not\in L$, then $k_1!\cdots k_t!Q_w$ has trivial image in $\Sh(X)_{\indec,s}$, and therefore $Q_w$ has trivial image in $\Sh(X)_{\indec,s}\tensor(\dbZ/r)$.

Therefore the images of $\sum_{w\in X^s}\dbZ Q_w$ and $\sum_{w\in L}\dbZ Q_w$ in
$\Sh(X)_{\indec,s}\tensor(\dbZ/m)$ coincide.
But  by (\ref{Radford basis}), the former sum is $\sum_{w\in X^s}\dbZ w$, whereas by the construction of $Q_w$, the latter sum is $\sum_{w\in L}\dbZ w$.
Hence these two images are the full $\Sh(X)_{\indec,s}\tensor(\dbZ/r)$ and $\Sh(X)_{\indec,s}^{(L)}\tensor(\dbZ/r)$, respectively.
\end{proof}

We now take $\Phi=(\Lam\colon G\to\dbZ_p^{\times,1} ,\tau\colon L\to G,e)$ to be a $p$-adic Magnus formation, as in Example \ref{the shuffle Magnus formation}.
Thus $G$ is a free pro-$p$  group on basis $X$, and we recall that $\Phi$ is compatible with the infiltration product $\downarrow$.
We further assume that the map $e$ is given by $e(i)=p^{j(i)}$ for some map $j\colon\{1,2\nek n\}\to\dbZ_{\geq0}$.
The map $j$ is weakly decreasing and is not identically $0$, by Definition \ref{def of Magnus formations}(vii).
Thus $j(1)\geq1$.

\begin{thm}
\label{first isom for shuffle algebra}
Suppose that $n<p$ and $L$ is the set of all Lyndon words in $X^*$ with length in $I_e$.
Then there is a canonical isomorphism of $\dbZ/p$-linear spaces
\[
\bigl(\bigoplus_{s\in I_e}\Sh(X)_{\indec,s}\bigr)\tensor(\dbZ/p)\xrightarrow{\ \sim\ }  H^2(G/G^\Phi).
\]
\end{thm}
\begin{proof}
By (\ref{indec commutes with tensor product}), $\Sh(X)_\indec^{(L)}\tensor\dbZ_p\isom (\Sh(X)\tensor\dbZ_p)_\indec^{(L)}$.
Hence, by Proposition \ref{the tensorized shuffle module} (with $r=p$),
\[
\begin{split}
&\bigl(\bigoplus_{s\in I_e}\Sh(X)_{\indec,s}\bigr)\tensor(\dbZ/p)=\bigl(\bigoplus_{s\in I_e}\Sh(X)_{\indec,s}^{(L)}\bigr)\tensor(\dbZ/p) \\
=&\Sh(X)_\indec^{(L)}\tensor(\dbZ/p)
=(\Sh(X)\tensor\dbZ_p)_\indec^{(L)}\tensor(\dbZ/p).
\end{split}
\]
By Theorem \ref{indecomposable quotient as the space of maps} for $A=\Sh(X)\tensor\dbZ_p$ and $\downarrow$, the latter module is isomorphic to $\Hom(G^\Phi,\dbZ/p)^G$ via the map $\bar w\tensor1\mapsto\rho_w^0$, for $w\in L$.
Moreover, by the definition of the first cohomology group, $\Hom(G^\Phi,\dbZ/p)^G=H^1(G^\Phi)^G$.
We deduce that
\begin{equation}
\label{ttt}
\bigl(\bigoplus_{s\in I_e}\Sh(X)_{\indec,s}\bigr)\tensor(\dbZ/p)\isom H^1(G^\Phi)^G.
\end{equation}

Now when $w=(x)$ is a word of length $1$ we have $\sig_w=x^{p^{j(1)}}\in G^p$.
When $w$ is a Lyndon word of length $\geq2$ we have $\tau_w\in G^{(2)}$, by Definition \ref{def of Magnus formations}(v), whence also $\sig_w\in G^{(2)}$.
Thus $G^\Phi$ is contained in the Frattini subgroup $G^pG^{(2)}=G^p[G,G]$ of $G$.
It follows that the inflation map $H^1(G/G^\Phi)\to H^1(G)$ is an isomorphism.
Since $G$ is a free pro-$p$ group, $\cd_p(G)\leq1$.
The five term sequence in profinite cohomology \cite{NeukirchSchmidtWingberg}*{Prop.\ 1.6.7} therefore implies that the transgression map
\[
\trg\colon H^1(G^\Phi)^G\to H^2(G/G^\Phi)
\]
is an isomorphism, and we combine it with (\ref{ttt}).
\end{proof}

\begin{rem}
\label{remark on Massey products}
\rm
Explicitly, this isomorphism is given as follows:
The left-hand side is generated by elements of the form $\bar w\tensor1$, where $w$ is a Lyndon word of length $s\in I_e$, and such a generator is mapped to $\trg(\rho_w^0)$.
Let $\alp_w\in H^2(\dbU_w/\dbU_w^0)$ correspond to the central extension
\[
1\to\dbU_w^0(\isom\dbZ/p)\to\dbU_w\to\dbU_w/\dbU_w^0\to1
\]
under the Schreier correspondence  \cite{NeukirchSchmidtWingberg}*{Th.\ 1.2.4}.
Let $\bar\rho_w\colon G/G^\Phi\to\dbU_w/\dbU_w^0$ be the homomorphism induced by $\rho_{w,R_w}\colon G\to\dbU_w$ (see Proposition \ref{image of SL}(a)).
By \cite{Hoechsmann68}, $\trg(\rho_w^0)$ is the pullback $\bar\rho_w^*(\alp_w)$ of $\alp_w$ to $H^2(G/G^\Phi)$ along $\bar\rho_w$.
It corresponds to the central extension
\[
0\to\dbZ/p\to\dbU_w\times_{\dbU_w/\dbU_w^0}(G/G^\Phi)\to G/G^\Phi\to1,
\]
where the middle term is the fiber product.
\end{rem}

We examine these cohomology elements in two special situations:

\begin{exam}
\label{example with Bocksteins}
\rm
Suppose that $w=(x)$ has length $1$.
Then $\alp_w$ corresponds to the extension
\[
0\to\dbZ/p\to\dbZ/p^{j(1)+1}\to\dbZ/p^{j(1)}\to0.
\]
The \textsl{Bockstein map} $\Bock\colon H^1(G/G^\Phi,\dbZ/p^{j(1)})\to H^2(G/G^\Phi)$ is the connecting homomorphism associated to this short exact sequence of trivial $G/G^\Phi$-modules \cite{NeukirchSchmidtWingberg}*{Th.\ 1.3.2}.
The homomorphism $\bar\rho_w\colon G/G^\Phi\to\bar\dbU_w$ may be identified with the map $\bar\eps_{(x)}\eps_{(x)}\ (\hbox{mod } p^{j(1)})$.
Then the pullback $\bar\rho_w^*(\alp_w)$ is $\Bock(\bar\rho_w)=\Bock(\bar\eps_{(x)})$  \cite{Efrat17}*{Example 7.4(1)}.
\end{exam}

\begin{exam}
\label{example with Massey products}
\rm
Let $w=(a_1\cdots a_n)$ be a Lyndon word in $X^*$ of length $n$, and suppose that $j(n)=0$.
Then the pullbacks $\bar\rho_w^*(\alp_w)$ are elements of the $n$-fold Massey product
\[
\langle\cdot,\nek\cdot\rangle\colon H^1(G)^n\to H^2(G).
\]
We refer, e.g., to \cite{Efrat14} for the definition of this map in the context of profinite cohomology, and recall that this is a multi-valued map, i.e., $\langle\varphi_1\nek\varphi_n\rangle$ is a \textsl{subset} of $H^2(G)$.
Namely, in this case $R_w=\dbZ/p$, and $\bar \dbU_w$ is the group $\dbU_w$ with its $(1,n+1)$-entry deleted.
Let $\bar\Lam_w\colon G\to(\dbZ/p)\langle\langle X\rangle\rangle^{\times,1}$ be the homomorphism induced by $\Lam$ as in \S\ref{section on unitriangular representations},
and denote the coefficient of a word $u$ in $\Lam_w(\sig)$ by $\bar\eps_u(\sig)$.
As shown by Dwyer \cite{Dwyer75} in the discrete case (see \cite{Efrat14}*{Prop.\ 8.3} for the profinite case) the pullbacks $\bar\rho_w^*(\alp_w)$ are elements of  $\langle\bar\eps_{(a_1)}\nek\bar\eps_{(a_n)}\rangle$.
\end{exam}

Finally, we specify Theorem \ref{first isom for shuffle algebra} in the two special cases discussed in the Introduction:

\begin{exam}
\label{main result for lower p central filtration}
\rm
Consider the binomial map $e(i)=p^{n-i}$, as in Example \ref{first example}.
Then $I_e=\{1,2\nek n\}$, and $L$ contains all Lyndon words $w$ of length $1\leq i\leq n$.
For a free pro-$p$ group $G$ on the basis $X$, let $K^{(n,p)}$ be its closed subgroup generated by all powers $\tau_w^{p^{n-i}}$ for such $w$.
We obtain that when $n<p$,
\[
\bigoplus_{s=1}^n\Sh(X)_{\indec,s}\tensor(\dbZ/p)\xrightarrow{\sim}H^2(G/K^{(n,p)}).
\]

The groups $K^{(n,p)}$ are closely related to the lower $p$-central filtration $G^{(n,p)}$, $n=1,2\nek$ of $G$.
Namely, in  \cite{Efrat17}*{Th.\ 5.3} it is shown (using Lie algebra techniques) that  the subgroups $K^{(n,p)},G^{(n,p)}$ coincide modulo $G^{(n+1,p)}$.
\end{exam}

\begin{exam}
\label{main result for p Zassenhaus filtration}
\rm
For a positive integer $t$, we consider the binomial map $e(i)=p^{t\lceil\log_p(n/i)\rceil}$, as in Example \ref{second example}.
Assume that $n<p$.
By Proposition \ref{conditions for J(n)}, $I_e=\{1, n\}$.
In the free pro-$p$ group $G$ on basis the $X$, let $K_{(n,p)}$ be the closed subgroup generated by all powers $x^{e(1)}$, $x\in X$, and by all Lie elements $\tau_w$, where $w$ is a Lyndon word of length $n$ in $X^*$.
We obtain that
\[
(\bigoplus_{x\in X}(\dbZ/p))\oplus\bigl( \Sh(X)_{\indec,n}\tensor(\dbZ/p)\bigr)\xrightarrow{\sim}H^2(G/K_{(n,p)}).
\]

When $t=1$, the subgroups $K_{(n,p)}$ are closely related to the $p$-Zassenhaus filtration of $G$ (see the Introduction).
Namely, in  \cite{Efrat23}*{Th.\ 4.6} it is shown using $p$-restricted Lie algebra techniques that $K_{(n,p)},G_{(n,p)}$ coincide modulo $(G_{(n,p)})^p[G,G_{(n,p)}]$.
\end{exam}

\end{document}